\documentclass[11pt]{amsart}
\usepackage{amsmath}
\usepackage[left=1in,top=1in,right=1in,bottom=1in]{geometry}
\usepackage{color}
\usepackage{mathrsfs}
\usepackage{url}

\usepackage{graphicx} % Required for inserting images
\usepackage{amsmath,amsthm,amssymb}
\usepackage{cleveref}
\usepackage{tikz, tikz-cd}

\title{Representation stability for moduli spaces of admissible covers}
\author[M. Chang-Lee]{Megan Chang-Lee}
\address{Department of Mathematics, Brown University, Providence, RI, USA} 
\email{\url{megan_chang-lee@brown.edu}}
\author[S. Kannan]{Siddarth Kannan}
\address{Department of Mathematics, Massachusetts Institute of Technology, Cambridge, MA, USA}
\email{\url{spkannan@mit.edu}}
\author[P. Tosteson]{Philip Tosteson}
\address{Department of Mathematics, University of North Carolina, Chapel Hill, NC, USA}
\email{\url{ptoste@unc.edu}}
\newcommand{\tFSA}{\widetilde{\mathbf{FWS}}_A}
\newcommand{\FS}{{\mathbf{FS}}}
\newcommand{\tV}{\widetilde{V}}

\newcommand{\FSA}{{\mathbf{FWS}}_A}
\newcommand{\M}{{\mathcal M}}
\newcommand{\Vect}{\mathrm{Vect}}
\newcommand{\mc}[1]{\mathcal{#1}}

\newcommand{\glue}{\mathrm{ glue}}

\newcommand{\cV}{\mathcal V}
\newcommand{\op}{^{\rm op}}
\newcommand{\inv}{^{-1}}
\newtheorem{thm}{Theorem}[section]
\newtheorem{lemma}[thm]{Lemma}
\newtheorem{prop}[thm]{Proposition}

\newcommand{\Hom}{\operatorname{Hom}}
\newcommand{\Adm}{\mathrm{Adm}}
\newcommand{\Mbar}{\overline{\mathcal{M}}}
\newcommand{\cM}{\mathcal M}
\newcommand{\bbQ}{\mathbb Q}

\newcommand{\bbN}{\mathbb N}
\newcommand{\bbC}{\mathbb C}

\newtheorem{innercustomthm}{Theorem}
\newenvironment{customthm}[1]
  {\renewcommand\theinnercustomthm{#1}\innercustomthm}
  {\endinnercustomthm}
  
\newtheorem{innercustomcor}{Corollary}
\newenvironment{customcor}[1]
  {\renewcommand\theinnercustomcor{#1}\innercustomcor}
  {\endinnercustomcor}

\newcommand{\onto}{\twoheadrightarrow}
\newcommand{\oast}{\circledast}
\newcommand{\bC}{\mathbf{C}}
\theoremstyle{definition} % text is not in italics
\newtheorem{ex}[thm]{Example}
\newtheorem{remark}[thm]{Remark}
\newtheorem{defn}[thm]{Definition}
\newtheorem{notn}[thm]{Notation}

\DeclareMathOperator{\id}{id}
\DeclareMathOperator{\Spec}{Spec}

\newcommand{\mainbound}{g + 5i}
\begin{document}

\maketitle

\begin{abstract}
We prove a representation stability result for the sequence of spaces $\Mbar_{g, n}^A$ of pointed admissible $A$-covers of stable $n$-pointed genus-$g$ curves, for an abelian group $A$. For fixed genus $g$ and homology degree $i$, we give the sequence of rational homology groups $H_i(\Mbar_{g, n}^A;\bbQ)$ the structure of a module over a combinatorial category, \textit{\`a la} Sam--Snowden, and prove that this module is generated in degree at most $g + 5 i$. This implies that the generating function for the ranks of the homology groups is rational, with poles in the set $\left\{-1, -\frac{1}{2}, \ldots, -\frac{1}{|A|^2\cdot(g + 5i)}\right\}$. In the case where $A$ is the trivial group, our work significantly improves on previous representation stability results on the Deligne--Mumford compactification $\Mbar_{g, n}$. 
\end{abstract}

\section{Introduction}
    Let $A$ be a finite abelian group. We study the moduli space $\Mbar_{g, n}^A$ of pointed admissible $A$-covers, introduced in \cite{JKK}, from the perspective of representation stability. Recall that the moduli space $\Mbar_{g, n}^A$ parameterizes morphisms \[(P, \tilde{p_1}, \ldots, \tilde{p_n}) \to (C, p_1, \ldots, p_n) \]
    where $(C, p_1, \ldots, p_n)$ is a Deligne--Mumford stable $n$-pointed curve of genus $g$, and the morphism $P \to C$ is a principal $A$-bundle away from the marked points and nodes of $C$. The marked points $\tilde{p_i}$ on the source curve $P$ are lifts of the marked points $p_i$ on the target curve $C$. These data are subject to additional smoothability conditions at nodes; see Section \ref{sec:moduli-background} for a precise definition. Fixing the genus $g$, a homology degree $i$, and the group $A$, we are interested in the asymptotic behavior of the homology group $H_i(\Mbar_{g, n}^A; \bbQ)$ as $n$ goes to infinity. 

    \begin{customthm}{A}\label{main-upshot}
    Fix a finite abelian group $A$, genus $g \geq 0$, and homology degree $i \geq 0$. Assume that $(g, n) \neq (0,0)$. Then the generating function
\[ \mu_{i,g, A}(t) := \sum_{n} \dim_{\bbQ} H_i(\Mbar_{g, n}^A; \bbQ) \cdot t^n \] is rational, and has the form
\[ \mu_{i,g, A}(t) = \frac{P(t)}{\prod_{j = 1}^{(g + 5i)|A|^2} (1 + j t)^{d_j}} \]
for some polynomial $P(t)$ and nonnegative integers $d_j$.
    \end{customthm}

    \subsection{The category $\tFSA$} To prove Theorem \ref{main-upshot}, we construct an action of an appropriate combinatorial category on the homology groups of the moduli space. It is more natural to construct this action if we replace the parameter $n$ with a finite set $X$, so that the admissible cover is labelled by the elements of $X$. The monodromy of an admissible $A$-cover over each marked point in the target curve is specified by a function $\ell_X: X \to A$. We obtain an action of a category $\tFSA$ on the homology groups $H_i(\Mbar_{g, X}^A(\ell_X);\bbQ)$, where
    \[\Mbar_{g, X}^A(\ell_X) \subset \Mbar_{g, X}^A \]
    is the substack consisting of admissible covers with monodromy specified by $\ell_X$. The category $\tFSA$ is a generalization of the category $\FS$ of finite sets and surjections; indeed, when $A$ is the trivial group, we have $\tFSA = \FS$. The objects of $\tFSA$ are pairs $(X, \ell_X)$, where $X$ is a finite set, and $\ell_X: X \to A$ is any function. The morphisms are somewhat more complicated; at first approximation, they are surjections satisfying a compatibility condition, but with an additional twist. See Section \ref{sec:category-background} for a precise definition. The upshot of this construction is that the association
    \[(X, \ell_X) \mapsto H_i\left(\Mbar_{g, X}^A(\ell_X) ; \bbQ\right) \]
    can be viewed as a functor from $\tFSA\op$ to vector spaces, i.e. an $\tFSA\op$-module. The functorial structure comes from appropriate gluing maps on moduli spaces of admissible $A$-covers, as described in \cite[Section 2]{JKK} and explored from the operadic perspective in \cite{petersen-operad}.  Our main structural result is the following.
     \begin{customthm}{B}\label{main-technical}
    Fix $i, g \geq 0$. The association \[(X, \ell_X) \mapsto H_i\left(\Mbar_{g, X}^{A}(\ell_X) ; \bbQ\right)\]
    admits the structure of a finitely generated $\tFSA\op$-module, which is generated in degree at most $\mainbound$ if $(i, g) \neq (0, 0)$ and in degree $1$ if $(i, g) = (0,0)$.
    \end{customthm}

    We take the convention that $\overline \cM_{g,X}( \ell_X)$ is a point for $g = 0$ and $|X| =0,1,2$, and when $g = 1$ and $|X| = 0$, so that in these cases its rational homology is $\bbQ$, concentrated in homological degree $0$.   
   
    Theorem \ref{main-upshot} follows from Theorem \ref{main-technical}, together with Proposition \ref{mainrestrictionsummary} which relates $\tFSA \op$-modules to $\FS\op$-modules, and a result of Sam-Snowden on the Hilbert series of $\FS \op$-modules \cite[Corollary 8.1.4]{sam2017grobner}.  Theorem \ref{main-technical} has many additional consequences beyond Theorem \ref{main-upshot}.   For instance, it determines the structure of the sequence of vector spaces $\bigoplus_{\ell: [n] \to A} H_i\left(\Mbar_{g, n}^{A}(\ell) ; \bbQ\right)$   considered as a sequence of  $A^n$ representations by the theory of \cite{tosteson2024hilbert, sam2019representations} and as a sequence of $S_n$ representations by \cite{tosteson2024hilbert}. 

    An interesting case of Theorem \ref{main-technical} is when the group $A$ is trivial. In this case the moduli space $\Mbar_{g, X}^A = \Mbar_{g, X}$ recovers the Deligne--Mumford compactification of the moduli space of curves, and Theorem \ref{main-technical} specializes to the $\FS\op$-action on the homology of these moduli spaces, as studied by the third author in \cite{Tosteson2}. After this specialization, Theorem \ref{main-technical} is a significant improvement of the main result in \cite{Tosteson2}, where an $O(g^2i^2)$ bound on generation degree of the $\FS\op$-module $X \mapsto H_i(\Mbar_{g, X}; \bbQ)$ is given. This bound was improved to $O(gi + i^2)$ in \cite{TostesonStableMaps}. Our bound is linear, and this improvement comes from a novel conceptual argument for finite generation. While the bound on generation degree in \cite{Tosteson2} relied on a combinatorial analysis of the asymptotic behavior of stable graph strata in $\Mbar_{g, n}$, our bound comes from a simpler inductive argument.  The main inputs are Harer's theorem on the virtual cohomological dimension (vcd) of the mapping class group \cite{harer1986virtual} and the purity of Deligne's weight filtration on the rational cohomology of $\Mbar_{g, X}^A$.

\subsection{Unpointed admissible $A$-covers} Our results also allow us to deduce finite generation results for the moduli spaces $\overline{\Adm}_{g, n}^A$ of admissible $A$-covers: introduced by Abramovich--Corti--Vistoli \cite{ACV} and Abramovich--Vistoli \cite{AV}, these parameterize the same covers as $\Mbar_{g, n}^A$, except that there is no chosen lift of the marked points of the target curve. By taking coinvariants with respect to the $A^X$-action on $\Mbar_{g, X}^A(\ell_X)$, we obtain an $\FSA\op$-module structure on the homology groups $H_i(\overline{\Adm}_{g, X}^A(\ell_X); \bbQ)$. Here, $\FSA$ is a subcategory of $\tFSA$ with the same objects, defined precisely in Section \ref{sec:category-background}. It is a corollary of Theorem \ref{main-technical} that this $\FSA\op$-module is also finitely generated.
\begin{customcor}{C}
If $(i, g) \neq (0, 0)$, the $\FSA\op$-module
\[(X, \ell_X) \mapsto H_i(\overline{\Adm}_{g, X}^A(\ell_X);\bbQ) \]
is finitely generated in degree at most $|A|(g + 5i)$. If $(i, g) = (0, 0)$, it is finitely generated in degree $1$. 
\end{customcor}

To describe the numerical consequences of Corollary C, we associate a multivariate generating function to the space of admissible covers. Given a labeled set $(X, \ell)$, we associate to $X$ a monomial in the polynomial ring $\bbC[t_a, a \in A]$, with variables indexed by $a$, by defining $t^{X} = \prod_{x \in X} t_{\ell(x)}$.  This monomial determines the isomorphism class of the labeled set $X$.  Applying a result on the Hilbert series of finitely generated $\FSA\op$ modules (Theorem \ref{HilbertFWSA}) we immediately obtain.
\begin{customcor}{D}
    The generating function  $$\sum_{[X] \text{ isomorphism class of $A$-labeled set}} \dim(H_i(\overline\Adm_{g,X}^A(\ell_X), \bbQ)) ~t^X$$ is rational with denominator a product of terms of the form $$(1- \frac{j \zeta}{|A|} t_a)$$ where $a \in A$  and $\zeta$ is an ${\rm ord}(a)$th root of unity and $j \in \bbN, j \leq |A|^2$.  
\end{customcor}
Corollary D specializes to a single variable generating function statement for the generating function for $\dim(H_i(\overline\Adm_{g,X}^A;\bbQ))$, resembling Theorem A.   Theorem A also admits a multivariate generalization using Theorem \ref{HilbertFWSA}.

\subsection{Relation to other work}
Finite generation results for the homology of the moduli space of curves $\M_{g, n}$ as an $\mathbf{FI}$-module have been proven in \cite{JR-mgn}. Homological stability for $\M_{g, n}$ goes back to Harer \cite{harerstability}. Finite generation for $H_i(\Mbar_{g, n};\bbQ)$ as an $\FS\op$-module was proven by the third author in \cite{Tosteson2}, and a similar result was obtained for Kontsevich spaces of stable maps in \cite{TostesonStableMaps}. Pointed admissible $G$-covers were studied from the operadic point of view by Petersen \cite{petersen-operad} and were used to study the moduli space $\mathcal{H}_{g, n}$ of $n$-pointed hyperelliptic curves by Brandt--Chan--Kannan \cite{BCK}. The class of $\overline{\Adm}_{0, n}^{G}$ in the Grothendieck ring of stacks has been determined by Bagnarol--Perroni \cite{BagnarolPerroni}. Cycle classes on $\Mbar_{g, n}$ obtained from moduli spaces of admissible $G$-covers have been studied by Schmitt--van Zelm \cite{SchmittvanZelm}.

\subsection{Summary of paper}
We recall the definition of the moduli space $\Mbar_{g, n}^A$ and its monodromy stratification in Section \ref{sec:moduli-background}. In Section \ref{sec:category-background} we construct the combinatorial categories $\FSA$ and $\tFSA$, and collect some basic facts about Hilbert series of modules over these categories. In Section \ref{sec: restricted modules}, we characterize finite generation of modules over $\FSA$ and $\tFSA$ in terms of certain suspension operations. We then construct an $\tFSA\op$-module structure on the homology groups of $\Mbar_{g, n}^A$ in Section \ref{section: Actions}. Using an inductive argument and the criteria developed in Section \ref{sec: restricted modules}, we prove that these modules are finitely generated in Section \ref{sec:finite_generation}. We conclude the paper with some further questions in Section \ref{sec:further questions}.

\subsection*{Acknowledgments} We thank Melody Chan for many helpful conversations. SK is supported by an NSF Postdoctoral Fellowship (DMS-2401850).

\section{Pointed admissible $A$-covers}\label{sec:moduli-background}
Below we recall the definitions of the moduli space $\overline{\Adm}_{g,n}^A$ of admissible $A$-covers as introduced by Abramovich--Corti--Vistoli \cite{ACV}, and its pointed analogue $\Mbar_{g, n}^A$ introduced by Jarvis--Kimura--Kaufmann \cite{JKK}.  We work over $\mathbb{C}$ throughout.
\begin{defn}%Admissible cover
  Given an abelian group $A$, an \textit{$n$-marked admissible $A$-cover} is an $n$-marked, stable, genus-$g$ nodal curve $(C, p_1, \dots, p_n)$, along with a covering map of curves $\pi: E \to C$ such that 

    \begin{enumerate}
        \item $\pi$ maps nodes of $E$ to nodes of $C$,
        \item $\pi$ is a principal $A$-bundle away from nodes and markings of $C$, 
        \item At points $p \in P$ mapping to a node of $C$, the structure of the maps $P \to C$ is given locally by: 
        \[\Spec \mathbb{C}[z,w]/(zw) \to \Spec \mathbb{C}[x,y]/(xy)\]
        where $x = z^r$ and $y = w^r$, 
       
        \item For each node $q \in P$, the action of the stabilizer of $q$, $A_q \subset A$, is \textit{balanced}: for each $a \in A_q$ the eigenvalues of the action of $a$ on the two tangent spaces at $q$ are inverses. 
    \end{enumerate}
\end{defn}

Conditions (3) and (4) guarantee the smoothability of the cover. Admissible $A$-covers of fixed genus $g$ and number of marked points $n$ form a smooth Deligne Mumford stack, which we will denote $\overline{\Adm}_{g,n}^A$. %thm 3.0.2 plus 4.3.2 in ACV

\begin{defn}%mgn^A bar 
    A \textit{pointed $n$-marked admissible $A$-cover} \[\pi: (E, \tilde{p}_1, \dots, \tilde{p}_n)\to (C, p_1, \dots, p_n) \] is an $n$-marked admissible $A$-cover $(\pi: E \to C, p_1, \dots, p_n)$ along with a choice of $n$ marked points $\tilde{p}_i \in E$ such that $\pi(\tilde{p}_i) = p_i$ for all $i$. 
\end{defn}

Pointed admissible $A$-covers of genus $g$ with $n$ marked points also form a stack, which we denote $\overline{\mc{M}}_{g,n}^A$. 

\begin{remark}
      Given a pointed admissible $A$-cover  \[\pi: (E, \tilde{p}_1, \dots, \tilde{p}_n)\to (C, p_1, \dots, p_n),\]
    there is a well-defined monodromy element $m_i \in A$ associated to each lifted point $\tilde{p}_i$. Given a loop $\gamma_i \subset C \smallsetminus \{p_1, \dots, p_n\}$ around $p_i$, pick a lift $\tilde{\gamma_i}$ of $\gamma_i$ which is in a small neighborhood of $\tilde{p_i}$, say starting at $q \in E$. Then $\tilde{\gamma_i}$ will end at a point $m_i \cdot q$, where $m_i \in A$ is a fixed group element, independent of the choice of $\tilde{\gamma_i}$. This element $m_i$ is the monodromy around $p_i$. Changing the choice of $\tilde{p}_i$ to $g \cdot \tilde{p}_i$ changes $m_i$ by conjugation by $g$. Since $A$ is assumed to be abelian, we see that in fact $m_i$ does not depend on the choice of lift $\tilde{p}_i$.
\end{remark}

\begin{defn}%particular monodromy
    For $(m_1, \dots, m_n) \in A^n$, let  
    \[\overline{\Adm}_{g,n}^A(m_1, \dots, m_n) \subset \overline{\Adm}_{g,n}^A \]
    denote the substack of $n$-marked admissible $A$-covers with monodromy element $m_i$ over $p_i$, and let
    \[\overline{\mc{M}}_{g,n}^A(m_1, \dots, m_n) \subset \Mbar_{g,n}^A \] denote its pointed analogue. 
\end{defn}

By Theorem 2.4 of \cite{JKK}, both $\overline{\mc{M}}_{g,n}^A$ and $\overline{\mc{M}}_{g,n}^A(m_1, \dots, m_n)$ are smooth Deligne-Mumford stacks. 

\begin{notn}
It is straightforward to replace the marking set $\{1, \ldots, n\}$ in the definition of $\overline{\Adm}_{g,n}^A$ or $\Mbar_{g,n}^A$ with an arbitrary finite set $X$. We use the notation $\overline{\Adm}_{g,X}^A$ and $\Mbar_{g, X}^A$ for the resulting moduli spaces. Similarly, the tuple $(m_1, \ldots, m_n)$ of monodromy values may be replaced with a function $\ell_X : X \to A$. We use the notation $\overline{\Adm}_{g,X}^A(\ell_X) \subset \overline{\Adm}_{g,X}^A$ and  $\Mbar_{g, X}^{A}(\ell_X) \subset \Mbar_{g,X}^A$ for the corresponding strata.
\end{notn}

Ultimately we will view the collection of groups $\{H_i(\Mbar_{g,X}^A(\ell_X);\bbQ) \}$ as a module over a combinatorial category, as $X$ varies over all finite sets and $\ell_X : X \to A$ varies over all possible functions.

\subsection{Rational cohomological dimension of moduli spaces}
    A key input to our inductive argument for finite generation will be Harer's theorem \cite{harer1986virtual} on the vcd of of the mapping class group. Recall that a group $G$ (or a space $X$) is said to have rational cohomological dimension at most $d$ if for every rational representation $V$ of $G$  (respectively, rational local system $V$ on $X$), the cohomology $H^i(G,V)$ (respectively $H^i(X,V)$) vanishes for all $i > d$. If $X$ is aspherical (i.e. is the classifying space of its fundamental group), then the cohomological dimension of $X$ equals the cohomological dimension of $\pi_1(X)$.

    \begin{thm}[Harer \cite{harer1986virtual}] \label{Harer}
        The orbifold fundamental group of $\cM_{g,n}$ has rational cohomological dimension 
        \[  \begin{cases} n-3 & g = 0, n \geq 3 \\ 4g-5 & g \geq 1,  n = 0 \\ 4g+n-4 & g \geq 1, n \geq 1 \end{cases} \]
    \end{thm}
Let $H_*^{\mathrm{BM}}$ denote Borel--Moore homology. We use Harer's theorem and Poincar\'e duality to denote the following vanishing result.
    \begin{prop}\label{vanish}
        Let $f: X \to \cM_{g,n}$ be a finite \'etale cover.  Then $H_i^{\rm BM}(X, \bbQ)$ vanishes for all $i \leq$
        \[\begin{cases} 2g +n - 3 & \text{ if } g \geq 1, n \geq 1 \\
                           n-4    & \text{ if } g = 0, n \geq 3 \\
                          2g - 1  & \text{ if } g \geq 2, n = 0 \\

        \end{cases} \]
    \end{prop}
    \begin{proof}
        Suppose that $g \geq 2$ and $n \geq 1$.  The rational cohomology of $X$ equals the cohomology of the fundamental group of $\cM_{g,n}$ with coefficients in the representation  corresponding to the local system $f_* \bbQ_X$ on $\cM_{g,n}$.  Here, we consider $\cM_{g,n}$ as a stack, so the fundamental group is the mapping class group ${\rm Mod}_{g,n}$  and $f_* \bbQ_X$ corresponds to a representation $V$ of the mapping class group, because $\cM_{g,n}$ is a classifying space for its fundamental group. So by Harer's  result,  Theorem \ref{Harer}, we obtain that $H^j(X, \bbQ) = H^j({\rm Mod}_{g,n}, V) = 0$ for all $j \geq 4g + n - 3$.  The real dimension of $X$ is $6g-6 + 2n$, so by Poincar\'e duality we have that $H_i^{BM}(X, \bbQ) = H^{6g - 6 + 2n - i}(X, \bbQ)$  vanishes for all $i$ such that $6g-6 + 2n - i \geq 4g + n - 3$. The other cases are obtained similarly.
    \end{proof}

 \begin{remark}
    An alternate justification of Proposition \ref{vanish}, that avoids appealing to the stack theoretic fundamental group of $X$ is as follows.  Write $\cM_{g,n}$ as the stack theoretic quotient of the scheme $\cM_{g,n}[3]$ parameterizing marked curves with level $3$ structure by $\Gamma:=Sp_{2g}(\mathbb Z/3)$. Then $X$ is the quotient of $X[3]:= X \times_{\cM_{g,n}} \cM_{g,n}[3]$ by $\Gamma$.  Therefore $H^*(X, \bbQ) = H^*(X[3],\bbQ)^\Gamma.$  Then as a finite cover of $\cM_{g,n}[3]$ we have that $X[3]$ is a (disjoint union of) classifying space(s) for finite index subgroups of the level three mapping class group ${\rm Mod}_{g,n}[3] \subseteq {\rm Mod}_{g,n}$.   Any finite index subgroup of a group of virtual cohomological dimension $c$ also has virtual cohomological dimension $c$.  Thus we obtain the same vanshing of $H^j(X, \bbQ)$.
 \end{remark}

\section{Background on representations of categories}\label{sec:category-background}
We now introduce the combinatorial categories which will act on the homology of the moduli spaces $\overline{\Adm}_{g,n}^A$ and $\Mbar_{g,n}^A$.
\subsection{The categories $\FSA$ and $\tFSA$}
Let $A$ be an abelian group. For us there will be two relevant categories associated to $A$, denoted by $\FSA$ and $\tFSA$. In both cases, the objects of the underlying category are pairs $(X, \ell_X)$ where $X$ is a finite set and
\[ \ell_X: X \to A \]
is a function, called an $A$-\textit{labeling} of $X$. These categories are defined so that the map

\begin{equation}\label{eq: Vig}
\mathcal{V}_{i,g}: (X, \ell_X) \mapsto H_i\left(\Mbar_{g, X}^A(\ell_X);\bbQ\right)
\end{equation}
can be viewed as a functor 
\begin{equation}\label{eq:moduleasfunctor}
\tFSA^{\mathrm{op}} \to \Vect_{\bbQ}.
\end{equation} As is common in the literature on representation stability \cite{sam2017grobner}, we refer to a functor as in (\ref{eq:moduleasfunctor}) as an $\tFSA\op$-module. Given an $\tFSA\op$-module $M$, there is a natural way to produce a module of coinvariants $M_A$ such that
\[M_A(X, \ell_X) := M(X, \ell_X)_{A^X}. \]
The map $M_A$ is viewed as a functor
\begin{equation}\label{eq:fSasfunctor}
\FSA\op \to \Vect_{\bbQ},
\end{equation}
so we call it an $\FSA\op$-module. In the case of the module $\mathcal{V}_{i, g}$ of (\ref{eq: Vig}), taking coinvariants yields 
\[(\mathcal{V}_{i, g})_A(X, \ell_X) = H_{i}\left(\overline{\Adm}_{g, X}^A(\ell_X);\bbQ\right). \]
\begin{defn}
The category $\FSA$ is defined follows: the objects are pairs $(X, \ell_X)$ where $X$ is a finite set, and $\ell_X: X \to A$ is a function. A morphism $f: (X, \ell_X) \to (Y, \ell_Y)$ is a surjection $f: X \to Y$ such that for all $y \in Y$, we have
\[ \sum_{x \in f^{-1}(y)} \ell_X(x) = \ell_Y(y). \]
\end{defn}

\begin{defn}
The category $\tFSA$ is defined follows: the objects are the same as those of $\FSA$. Now, a morphism $(X, \ell_X) \to (Y, \ell_Y)$ consists of a pair $(f, g)$ where $f: (X, \ell_X) \to (Y, \ell_Y)$ is a morphism in $\FSA$, and $g: X \to A$ is an arbitrary function, called a \textit{pointing}. The composition of morphisms  $(f_1,g_1): (X, \ell_X) \to (Y, \ell_Y)$ and $(f_2, g_2) : (Y, \ell_Y) \to (Z, \ell_Z)$ is defined by \[(f_2, g_2) \circ (f_1, g_1) = (f_2 \circ f_1, g_1 + g_2 \circ f_1).\]
It is straightforward to check that $\tFSA$ is a category.
\end{defn}

\begin{remark}
   An object of $\FSA$ or $\tFSA$ is a set $X$ equipped with a labeling $\ell_X$ of its elements by elements of $A$.  Occasionally for notational simplicity, we will suppress $\ell_X$ and simply refer to $X$ as an $A$-labelled set.
\end{remark}

\begin{defn}
In general, for a category $\mathbf{C}$, a $\mathbf{C}^{\mathrm{op}}$-module over a field $k$ is a contravariant functor $\mathbf{C} \to \Vect_k$.
\end{defn}

\begin{ex}
    For any object $c$ of $\bC$,  we write $P_c$ for the \emph{principal projective}  (or \emph{representable})  $\bC\op$-module,  defined by $P_c(d) = k \Hom_C(d,c)$.   By the Yoneda lemma, maps of $\bC\op$-modules from $P_{c} \to M$ correspond to elements $m \in M(c)$.  
\end{ex}

\begin{defn} %Finite Generation
    Given a $\bC\op$-module $M$ for $\bC = \FSA$ or $\tFSA$, we say that $M$ is \textit{(finitely) generated} in degree $\leq d$ if there is a (finite) list of objects $X_i \in \bC$ with $|X_i| \leq d$ and elements $m_i \in M(X_i)$ such that the elements $m_i$ generate $M$ under linear combinations and the action of transition maps.
    
\end{defn}

\begin{remark}
    If $M(X)$ is finite-dimensional for all $X$, then $M$ is finitely generated in degree $\leq d$ if and only if it is generated in degree $\leq d$.
\end{remark}

    \begin{remark}
      There is an equivalent characterization of generation in degree $\leq d$.   A module $M$ is generated in degree $\leq d$ if and only if for every $X \in \tFSA$ with $|X| > d$, the vector space $M(X)$ is spanned by elements of the form $(f,g)^*(w)$ for $(f,g): X \to Y$  a morphism in $\tFSA$ to a set $Y$ with $|Y| < |X|$ and $w \in M(Y)$. 
    \end{remark}

The Yoneda Lemma gives us the following equivalent definition of finite generation. This version of finite generation will be used in most of our proofs. 

\begin{prop}\label{prop:princ_proj}
A $\tFSA\op$-module $M$ is finitely generated in degree $\leq n$ if and only if there exists a surjection \[\bigoplus_{i = 1}^r P_{X_i} \to M\] where $|X_i| \leq n$ for all $i$.
\end{prop}

\begin{prop} \label{FSTilderestriction}
    If an $\tFSA \op$-module $M$ is generated in degree $\leq d$, then its restriction to $\FSA \op$ is generated in degree $\leq d |A|$.  Conversely, if $M|_{\FSA\op}$ is generated in degree $\leq d$, then $M$ is generated in degree $\leq d$.
\end{prop}
\begin{proof}
    The converse statement is immediate, because if a list of elements generates the restriction $M|_{\FSA\op}$ as an $\FSA\op$-module then they generate $M$ as an $\tFSA \op$-module.
    
    For the forward direction, by passing to quotients it suffices to prove that for $X \in \tFSA$ with $|X| \leq d$, the restriction of the $\tFSA\op$-module $P_{X}$ to $\FSA\op$ is generated in degree $\leq d|A|$. We have that $P_X(Y) = \bbQ \tFSA(Y,X).$ Now a map  $Y \to X \in \tFSA$ corresponds to a map of labeled sets $(f,g): Y \to A \times X$  where $\ell_{A \times X} = \ell_X \circ \pi_X$, such that $\pi \circ (f,g)$ is a surjection of $\FSA \op$-modules $$\bigoplus_{S \subseteq A \times X, ~ \pi_X(S) = X} P_S \to P_X|_{\FSA\op},$$ taking the generator of $P_S$ to the canonical inclusion $S \to A \times X$.
\end{proof}

Let $\FS_A \subseteq \tFSA$ denote the subcategory spanned by $(X, 0)$, where $0: X \to A$ is the zero labeling. $\FS_A$ may be thought of as an unlabeled version of $\tFSA$. There is a functor $u: \FSA \to \FS_A$, given by $(X, \ell) \mapsto (X,0)$ and $(f,g) \mapsto (f,g)$.

\begin{remark}
    The category $\FS_A$ is equivalent to the category of finite free $A$-sets and surjections between them. $\FS_A$ is also known as the category of sets and $A$-surjections, see \cite{sam2019representations}.  
\end{remark}

Given an $\FSA \op$-module $M$, we define $u_*(M)$ to be the $\FS \op$-module $u_*(M)(X) := \prod_{\ell: X \to A} M(X,\ell)$.  A surjection $f: X \to Y$ acts by taking the factor of $u_*(M)(Y)$ corresponding to $\ell: Y \to A$ to the product $$\prod_{\tilde \ell : X \to A ~|~\left(\sum_{x \in f\inv(y))} \tilde \ell(x)\right) = \ell(y)} M(X, \tilde \ell)$$ via $\prod_{\tilde \ell} f^*$, where for each $\tilde \ell$ the map $f$ is considered as a map of labelled sets from $(X,\tilde \ell) \to (Y, \ell)$.

We have a commutative diagram of functors
\begin{center}
\begin{tikzcd}
\tFSA \arrow[r, "u"]               & \mathbf{FS}_A      \\
\FSA \arrow[r, "u"] \arrow[u, "i"] & \FS \arrow[u, "i"],
\end{tikzcd}
\end{center}
where $u$ is the functor that forgets the labeling of a set, and $i$ is the functor that takes a morphism $f$ to the pair $(f,0)$. The construction $u_*$ can be interpreted as the pushforward (or right Kan extension) of an $\FSA \op$ module along $u$.

\begin{remark}
    There are also functors going the other direction: 
        \begin{center}
        \begin{tikzcd}
        \tFSA \arrow[d, "p"] & \mathbf{FS}_A \arrow[l, "v"] \arrow[d, "p"] \\
        \FSA                 & \FS \arrow[l, "v"] ,                
        \end{tikzcd}
        \end{center}
    where $v$ labels every set by $0$ and $p$ takes a morphism specified by the pair $(f,g)$ to $f$.  These morphisms satisfy $u \circ v = \id, p \circ i = \id$.  (We will not make signficant use of these functors). 
\end{remark}

In fact, the functors $i^*: {\rm Mod}(\FS_A \op) \to {\rm Mod}(\FS\op)$ and $u_*: {\rm Mod}(\FSA\op) \to {\rm Mod}(\FS\op)$ are related.  We have the following Fourier equivalence, observed by Sam--Snowden \cite{sam2019representations}.

\begin{prop}[Fourier duality] \label{Fourier}
    Suppose that $|A|$ is invertible in $k$ and $k$ contains the $|A|$th roots of unity.   Let $A^{\vee} = \Hom(A, \mu)$ be the group of characters of $A$.  Then there is an equivalence of categories $$\mathbb F: {\rm Mod}(\FS_A \op) \to {\rm Mod}(\mathbf{FWS}_{A^{\vee}} \op), $$ from $\FS_A \op$-modules to $\mathbf{FWS}_{A^\vee} \op$-modules defined defined by $$\mathbb F(M)(X, \ell) = \Hom_{A^X}(k_\ell, M(X)),$$ where $k_{\ell}$ is the one-dimensional $A^X$ representation $\otimes_{x \in X} k_{\ell(x)}$  given by tensoring together the characters $\ell(x)$. 
    Furthermore we have that $u_* \circ \mathbb F$ is naturally isomorphic to $i^*$, and $M$ is generated in degree $\leq d$ if and only if $\mathbb F(M)$ is generated in degree $\leq d$.
\end{prop}
\begin{proof}
    For the equivalence of categories see  \cite[Lemma 6.6.4]{sam2019representations}.  The point is that an action by $A^X$ on vector space $V$ corresponds to a grading of $V$ by isotypic components $ \ell \in ( A^{\vee})^X$.  By definition we have that $$u_* \circ \mathbb F (M)(X)  = \bigoplus_{\ell \in ( A^{\vee})^X} \Hom_{A^X}(k_\ell, M(X)) \cong M(X) = i^*(M)(X)$$ and the isomorphism is compatible with the action of surjections $X' \to X$.
    
    For the statement on generation degree, we observe that if $X \in \FS_A$ then for the principal projective $P_X$ we have that its Fourier transform is a sum of principal projectives $$\mathbb F(P^{\FS_A\op }_X) = \bigoplus_{\ell : X \to A^\vee } P^{\mathbf{FWS}_{A^\vee}\op }_{(X, \ell)}.$$ The statement now follows from Proposition \ref{prop:princ_proj}. 
\end{proof}
 
\begin{prop} \label{restricting_FWS} Suppose that $k$  has characteristic relatively prime to $|A|$ and that $M$ is an $\FSA\op$-module.
If $M$ is generated in degree $\leq d$, then $u_* M$ is generated in degree $\leq |A|d$.
\end{prop}
\begin{proof} We may extend scalars to the algebraic closure of $k$ without affecting finite generation.   We have that $u_*(M) = i^*   \mathbb F^{\inv}(M).$  Hence by Proposition  \ref{Fourier}, it suffices to prove that the restriction of an $\FS_A\op$-module generated in degree $\leq d$ to $\FS$ is generated in degree $\leq d |A|$.   This is proved identically to Proposition \ref{FSTilderestriction}:  we have that $i^*(P_X)$ admits a surjection from $$\bigoplus_{S \subseteq X \times A~|~ \pi_X(S) = X} P_{X\times A},$$
so Proposition \ref{prop:princ_proj} gives the desired statement. 
\end{proof}

Combining Propositions \ref{FSTilderestriction} and \ref{restricting_FWS} we obtain the following relationship between $\tFSA \op$-modules and $\FS \op$-modules.

\begin{prop}\label{mainrestrictionsummary}
    If $M$ is an $\tFSA\op$-module generated in degree $\leq d$, then $X \mapsto \bigoplus_{\ell: X \to A} M(X, \ell_X)$ is an $\FS \op$-module generated in degree $\leq d |A|^2$.
\end{prop}

Modules over the category $\FSA\op$  were studied by Sam and Snowden in \cite[Section 5]{sam2019representations}.  Using Gr\"obner theory, they proved that $\FSA \op$ is Noetherian, and that finitely generated representations have rational Hilbert series.

\begin{thm}[Sam--Snowden \cite {sam2019representations} Theorem 5.1.1]
    A submodule of a finitely generated $\FSA \op$-module is finitely generated.  
\end{thm}

\subsection{Convolution tensor products}
For any category $\bC$ together with a symmetric monoidal tensor product $\sqcup: \bC \times \bC \to \bC$, there is an associated Day convolution tensor product $\oast$ on the category ${\rm Rep}(\bC)$ of $\bC$ representations, defined as follows.  Given $\bC$ representations $M$ and $N$, we define $M \oast N$ to be the left Kan extension along $\sqcup$ of the $\bC \times \bC$ representation $(c_1, c_2) \mapsto M_{c_1} \otimes N_{c_2}$.    By definition, this tensor product is bilinear and right-exact.

Both the categories $\FSA \op$ and $\tFSA\op$ have a symmetric monoidal structure given by disjoint union:   $$(X,\ell_X) \sqcup (Y, \ell_Y) = (X \sqcup Y, \ell_X \sqcup \ell_Y).$$
Accordingly, there is an associated convolution tensor product of $\FSA\op$ (resp. $\tFSA\op$)-modules.  In this case, the convolution tensor product can be defined concretely as

$$(M \oast N)(X,\ell_X) = \bigoplus_{X = A \sqcup B} M(A, \ell_X|_A) \otimes  N(B, \ell_X|_B),$$

where the sum is over all decompositions of $X$ into two disjoint subsets.  A surjection of labeled sets $f:Y \to X$ acts on $M \oast N$ by taking the summand corresponding to $A \sqcup B$ to the summand corresponding to $f\inv(A) \sqcup f\inv(B)$ via the map $$(f|_{f\inv(A)})^* \otimes (f|_{f\inv(B)})^* : M(A) \otimes N(B) \to M(f\inv(A)) \otimes N(f\inv(B)).$$

\begin{prop}\label{product}
    Let $M$ and $N$ be $\tFSA\op$-modules. If $M$ is finitely generated in degree $\leq d_M$ and $N$ is finitely generated in degree $\leq d_N$, then $N \oast M$ is finitely generated in degree $\leq d_N + d_M$.
\end{prop}

\begin{proof}
    Given labeled sets $X$ and $Y$, we have that $P_X \oast P_Y = P_{X \sqcup Y}$.  A $\tFSA\op$-module $M$ is  finitely generated in degree $\leq d_M$ if and only if there is a finite list $X_1, \dots, X_k$ of $A$-labelled sets such that $|X_k| < d_M$, and a surjection $\bigoplus_{i = 1}^k P_{X_i} \onto M$. Tensoring this surjection with the analogous surjection $\bigoplus_{i = 1}^\ell P_{Y_j} \onto N$,  we obtain a surjection $$\bigoplus_{(i,j) = (1,1)}^{(k,\ell)} P_{X_i \sqcup Y_j} \onto M \oast N, $$ using right exactness and bilinearity of the convolution tensor product.
\end{proof}
The same proposition holds and the same proof works for $\FSA\op$-modules as well.

\subsection{Hilbert Series}\label{subsec:HilbertSeries}
Let $R_A := \bbQ[[t_a, a \in A]]$ be the ring of formal power series in the variables $t_a, a \in A$. Given a function $f: A \to \bbN$ we let $t^f \in R_A$ be the monomial $\prod_{a \in A} t_a^{f(a)}$.  We write $C_f$ for the multinomial coefficient $\frac{|f|!}{\prod_{a \in A} f(a)!}$.

Let $M$ be an $\FSA\op$-module.  Following Sam--Snowden we associate to $M$ a generating function $H_M \in R_A$ that records the dimensions of the vector spaces $M(X, \ell_X)$, where $(X, \ell_X)$ ranges over labeled finite sets.  Given a function $f: A \to \bbN$ we write $X_f$ for the $\sqcup_{a \in A} [f(a)]$,  equipped with its canonical $A$-labeling $\ell_f$ defined so that $\ell_f([f(a)]) = \{a\}$.  As $f: A \to \bbN$ varies, $(X_f, \ell_f)$ ranges over all isomorphism classes of objects of $\FS_A$.   With this motivation, we define
\[ H_M := \sum_{f: A \to \bbN} \dim M(X_f, \ell_f)~  C_f ~ t^f.  \]
The following theorem was proved by Sam--Snowden \cite{sam2019representations}.
\begin{thm}\label{FWSHilbertSeries}
    For any finitely generated $\FSA\op$-module $M$, the Hilbert series $H_M$ is rational with denominator equal to a product of terms of the form $(1 + \sum_{a \in A} c_a t_a)$ where $c_a \in \bbQ(\zeta_{|A|})$ is a linear combination of roots of unity of order dividing $|A|$.
\end{thm}

To obtain Corollary D from Theorem C, we essentially apply Theorem \ref{FWSHilbertSeries}.  Note however, that the generating function of Corollary D is defined differently from $H_M$.  Following \cite{tosteson2024hilbert}, we describe an different generating function associated to an $\FSA\op$ module and the form it takes.    We let $\tilde H_M = \sum_{f: A \to \bbN} \dim M(X_f, \ell_f) t^f$.   (In fact there is another natural generating function $E_M = \sum_{f: A \to \bbN} \dim M(X_f, \ell_f) C_f \frac{t^f}{|f|!}$ that we will not use).  These generating functions determine each other, and from this it is straightforward to see that Theorem \ref{FWSHilbertSeries} implies that $\tilde H_M$ is rational when $M$ is finitely generated.  We have the following effective refinement.

\begin{thm}\label{HilbertFWSA}
  Let $M$ be an $\FSA\op$-module that is finitely generated in degree $\leq d$. Then the Hilbert series $\tilde H_M$ is rational with denominator equal to a product of terms of the form $(1 + \frac{j \eta}{|A|} t_a)$ where $\eta$ is an ${\rm ord}(a)$th root of unity, $j \in \bbN$, and $j \leq |A|^2$.
\end{thm}
\begin{proof}
   This follows from \cite[Theorem 2.1]{tosteson2024hilbert} applied to the Fourier dual $\FS_{A^{\vee}}\op$ module to $M$, using Proposition \ref{Fourier}.
\end{proof}

\section{Shifted modules and finite generation}\label{sec: restricted modules}

In this section we will give a technical reformulation of finite generation (Proposition \ref{SuspensionMapProp}) that we will use in order to prove the main theorem. To do so, we use the disjoint union of labelled sets to define shifted $\FSA\op$- and $\tFSA\op$-modules.  
\begin{defn}
       Let $M$ be a $\bC\op$-module, for $\bC= \FSA$ or $\tFSA$.  Given $a \in A$, we define the \textit{shift of} $M$ on a set $X$ with labeling $\ell_X$ by \[\Sigma_a M(X, \ell_X) := M((*, a) \sqcup (X, \ell_X)),\] where $(*,a)$ denotes the one-object set labelled by $a$.  

       This assignment also defines a functor from the category of $\tFSA\op$-modules to itself. It acts on morphisms of $\tFSA\op$-modules $(f,g)$ as follows: for $M$ an $\tFSA\op$-module and sets $(X, \ell_X)$ and $(Y, \ell_Y)$ in $\tFSA\op$, $\Sigma_a f: X \sqcup \{*\} \to Y \sqcup \{\bullet\}$ maps $X$ to $Y$ via $f$ and sends $* \mapsto \bullet$. Then $\Sigma_a g: X \sqcup \{*\} \to A$ send $X$ to $A$ via $g$ and sends $* \mapsto a$. 
\end{defn}

\begin{prop}\label{restricted}
     Let $M$ be a $\tFSA \op$-module. If $M$ is finitely generated in degree $\leq d$, then $\Sigma_a M$ is finitely generated in degree $\leq d$.
\end{prop}
\begin{proof}
     The process of shifting an $\tFSA \op$-module defines an exact functor $\Sigma_a$. By exactness, it suffices to prove the proposition in the case where $M = P_X$ for $(X, \ell_X) \in \tFSA$ a set of size $\leq d$.  For any $(Y, \ell_Y) \in \tFSA$, we have that $\Sigma_aP_X(Y)$ is the free abelian group on the set of $\tFSA$ morphisms  $(Y \sqcup{*}, \ell_Y \sqcup a) \to (X, \ell_X)$.

     Given an object $(X, \ell_X)$ in $\tFSA$ and fixing group elements $a, b \in A$, we construct the following elements of $\Sigma_a P_X$, which we claim generate $\Sigma_aP_X$ as a module. 
     \begin{enumerate}
        \item  For each $x \in X$ such that $\ell_X(x) = a$, we may consider the element $\id|_{X - x} \sqcup x:(X - x \sqcup *, \ell_X \sqcup a) \to (X, \ell_X)$ of $\Sigma_a P_X(X - x, \ell_{X}|_{X - x})$. As an $\tFSA$ morphism, this map sends $*$ to $x$ and $X - x$ to itself, and the pointing is given by $0 \sqcup b: X - x \sqcup * \to A$, which sends $*$ to $b$ and all elements of $X - x$ to $0$.  

        \item For each $x \in X$, we may consider the object $(X,\tilde \ell_X)$ of $\tFSA$ where $\tilde \ell_X(x) = \ell_X(x) - a$ and $\tilde \ell_X(y) = \ell_X(y)$ if $y \neq x$. There is an element of $\Sigma_aP_X(X, \tilde \ell_X)$ given by the $\tFSA$ morphism $\id_X \sqcup x : X \sqcup * \to X$  and the map $0 \sqcup b : X-x  \sqcup * \to A$.
     \end{enumerate}
     There are $\leq 2|X||A|$ such elements, and each has degree $\leq d$.  To prove that they generate $\Sigma_a P_X(Y)$, let $(Y \sqcup{*}, \ell_Y \sqcup a) \to (X, \ell_X)$ be an arbitrary morphism in $\tFSA$, represented by a pair $(f,g)$ where $f$ is a morphism $(Y, \ell_Y) \to (X, \ell_X)$ in $\FSA$ and $g: Y \to A$ is a function.

     It suffices to show that $(f,g)$ factors as a composite of a morphism $(\overline f \sqcup *, \overline g \sqcup 0)$ with one of two types of elements of $\Sigma_a P_X(Y)$ described above. Since $f$ is a surjection of sets, ${\rm im}(f|_Y)$ is either all of $X$ or ${\rm im}(f|_Y) = X - f(*)$.  In the first case, we factor $(f,g)$ as $$(f|_Y \sqcup *, g|_Y \sqcup 0): (Y \sqcup *, \ell_Y \sqcup a) \to (X \sqcup *, \tilde \ell_X \sqcup a)$$ followed by the element of type (2) associated to $x = f(*)$ and $b = g(*)$.     In the second case, we factor $(f,g)$ as $$(f|_Y \sqcup *, g|_Y \sqcup 0): (Y \sqcup *, \ell_Y \sqcup a) \to (X - f(*) \sqcup *,  \ell_X \sqcup a)$$ followed by the element of type (1) associated to $x = f(*)$ and $b = g(*)$.
     \end{proof}
 
\begin{prop}\label{quotient}
    Let $f: M \to N$ be a map of $\bC\op$-modules.  If $M$ is generated in degree $\leq d$ and $M/N$ is   generated in degree $\leq e$, then $N$ is generated in degree $\leq \max(d,e)$.
\end{prop}
\begin{proof}
     Let $\{n_i\}$ and $\{q_j\}$ be elements that generate $N$ and $M/N$ respectively. Then choosing lifts $\tilde q_j$ of $q_j$ to $M$, we have that $\{f(n_i)\} \cup \{\tilde q_j\}$ generate $M$.  So if $n_i$ are of degree $\leq d$ and $q_j$ are of degree $\leq e$, we have constructed generators of $M$ of degree $\leq \max(d,e)$.
\end{proof}

We define a $\tFSA\op$-module $\tV_0$ by 

\begin{equation}\label{eqn:tV0}
    \tV_0(X) =  \begin{cases}  \bbQ \{A^X/A\} & |X| \geq 3 \text{ and $\sum_{x \in X} \ell_X(x) = 0$} \\ 0 & \text{otherwise}, \end{cases}
\end{equation}
where a map $(f,g): (Y,\ell_Y) \to (X, \ell_X)$ acts via (the linearization of the quotient of) the map 
\[ A^X \to A^Y \qquad (h: X \to A) \mapsto ([h \circ f + g] : Y \to A) \]
if $\sum_{x \in X} \ell_X(x) = 0$ and acts by $0$ otherwise.

In other words, thinking of $A^X$ and $A^Y$ as the sets of $A$-labelings of $X$ and $Y$, respectively, we label the elements of $Y$ by their labels under the image of $f$, and then twist this labeling by $g$.

The $\tFSA\op$-module $\tV_0$ is closely related to one coming from the geometry of moduli spaces of pointed admissible covers. In particular, Proposition 2.4 of \cite{BCK} gives a classification of connected components of $\overline{\M}_{0,X}^A (\ell_X)$. As such, we define the $\tFSA\op$-module $\overline{V}_0(X)$ as, 
\begin{equation}\label{eqn:v0bar}
    \overline V_0(X) =  \begin{cases}  \bbQ \left(\frac{A}{\langle \ell_X(x)\mid  x \in X\rangle}\right)^X/A & \text{if }|X| \geq 3 \text{ and $\sum_{x \in X} \ell_X(x) = 0$} \\ 0 & \text{otherwise.} \end{cases}.
\end{equation}
Note that $\overline{V}_0$ receives a quotient map from $\tV_0$. For any $\tFSA\op$-module $M$, we define a map of $\tFSA \op$-modules

\begin{equation}\label{eta tilde}
\tilde \eta_M: \bigoplus_{a \in A} \Sigma_{-a} \tV_0 \oast \Sigma_a M \to M
\end{equation}
as follows.  We have that the left hand side evaluated on a set $(X, \ell_X)$ is \[\bigoplus_{a \in A} ~~ \bigoplus_{(X, \ell_X) = (S, \ell_S) \sqcup (T, \ell_T)} \tV_0(S \sqcup (*, -a)) \otimes M(T\sqcup (*,a))\]= \[  \bigoplus_{S \subseteq X, |S| \geq 2} \bbQ\{A^S\} \otimes  M\left((X-S, \ell_X|_{X-S}) \sqcup (*, \sum_{s \in S} \ell_X(s))\right).\] For every subset $S \subseteq X$ and every function $b \in A^S$ there is an associated morphism in $\tFSA$ \[(f_S, g_b): (X, \ell_X) \to \left((X-S) \sqcup *, \ell_X|_{X-S} \sqcup \sum_{s \in S}\ell_X(s) \right)\] 
where functions $f_S: X \to (X \setminus S)\sqcup \{*\}$ and $g_b: X \to A$ are given by
\begin{equation*}
    \begin{aligned}[c]
        f_S(x) = \begin{cases} * \text{ if } x \in S \\ x \text{ if } x \in X \setminus S \end{cases}
        \end{aligned}
        \qquad\qquad
        \begin{aligned}[c]
        g_b(x) = \begin{cases} b(x) \text{ if } x \in S\\ 0 \text{ if } x \in X \setminus S \end{cases}
    \end{aligned}
\end{equation*}

We define $\tilde \eta_M$ as the sum over all such maps where $|S| \geq 2$:
\[\tilde \eta_M =  \underset{S \subseteq X , |S| \geq 2,  b \in A^S} \sum (f_S, g_b)^*.\]

\begin{prop}\label{SuspensionMapProp}
    An $\tFSA\op$-module $M$ is generated in degree $\leq d$ if and only if $\tilde \eta_M$ is surjective when evaluated on $X$ for every $X \in \tFSA$ with $|X| > d$. 
    \end{prop}
\begin{proof}
     A module $M$ is generated in degree $\leq d$ if and only if for every object $X \in \tFSA$ with $|X| > d$, the vector space $M(X)$ is spanned by elements of the form $(f,g)^*(w)$ for $(f,g): X \to Y$  a morphism in $\tFSA$ to a set $Y$ with $|Y| < |X|$ and $w \in M(Y)$. 

     We may show that any morphism $(f,g): X \to Y$ for $|Y| < |X|$ factors through a morphism of the form $(f_S, g_b)$ as above for $|S| = 2$.  Indeed, let $s_1, s_2$ be two distinct elements mapped to the same element of $Y$, and let $S = \{s_1, s_2\}$.  Let $b: S \to A$ be defined by $b(s_i) = g(s_i)$. Then there is a unique factorization of $(f,g)$ as \[X \xrightarrow{(f_S, g_b)} (X- S) \sqcup *  \to Y. \] Consequently, $(f,g)^*(w)= (f_S, g_b)^*(w')$ where $w'$ is the pullback of $w$ by the action of the map $(X-S) \sqcup * \to Y$. 

     So an element of $M(X)$ is of the form $(f,g)^*(w)$ if and only if it is of the form $(f_S, g_b)^*(w')$ for some $S \subseteq X$ and $b: S \to A$, with $|S| \geq 2$. As mentioned above, elements of the form $(f,g)^*(w)$ span $M(X)$ if and only if $M$ is generated in degree $\leq d$.   On the other hand, by definition of $\tilde \eta_M$, elements of the form $(f_S, g_b)^*(w')$ span $M(X)$ if and only if $\tilde \eta_M$ is surjective.  
\end{proof}

\section{$\tFSA\op$-module and $\FSA\op$-module structures on homology groups}\label{section: Actions}
Fix nonnegative integers $i$ and $g$. We now describe how the assignment
\[\mathcal{V}_{i,g}: (X, \ell_X) \mapsto H_i(\Mbar_{g, X}^A(\ell_X);\bbQ) \]
of (\ref{eq: Vig}) defines an $\tFSA\op$-module when the space $\overline{\mc{M}}_{g, X}^A(\ell_X)$ is defined. When these spaces are not defined, when $|X| < 3$ and $g = 0$, we will define $\cV_{0,0}$ separately in Section \ref{sec:finite_generation} and set $\cV_{i,0} = 0$ if $i \neq 0$ and $|X| < 3$. 

The idea is as follows: given a morphism $(f, g): (X, \ell_X) \to (Y, \ell_Y)$ in $\tFSA$, we will define a map of spaces
\[ \varphi_{(f, g)^*} : \Mbar_{g, Y}^A(\ell_Y) \to \Mbar_{g, X}^A(\ell_X), \] by gluing a genus-$g$, $Y$-pointed admissible $A$-cover to a fixed genus-zero cover, which comes from the space
\[ \Mbar_{0, (f, g)}^A:= \prod_{y \in Y} \Mbar_{0, f^{-1}(y) \sqcup y}^A (\ell_X|_{f^{-1}(y)}, -\ell_Y(y)). \]
Though the map of spaces $\varphi_{(f, g)^*}$ depends on the choice of point from $\Mbar_{0, (f, g)}$, we will see that the connected component containing any choice of point is uniquely determined by $(f, g)$, so that we get a well-defined map on homology
\[ (f, g)^* :  H_*(\Mbar_{g, Y}^A(\ell_Y);\bbQ) \to H_*(\Mbar_{g, X}^A(\ell_X); \bbQ). \]
We then show that this system of maps has the structure of an $\tFSA\op$-module. We begin by reformulating the description of the connected components of the spaces $\Mbar_{0, X}^A(\ell_X)$ found in \cite{BCK}.

\begin{defn}
Given a point $(P \to C)$ of $\Mbar_{0, X}^A(\ell_X, a)$, define for each pair of points $x_1, x_2 \in X$ an element
\[ b_{x_1, x_2}(P \to C) \in A/ \langle \ell_X(x) \mid x \in X \rangle \]
as follows: let $P_{1}$ and $P_{2}$ denote the connected components of $P$ such that $\tilde{x_1} \in P_{1}$ and $\tilde{x_2} \in P_2$. Let $\tilde{b}_{x_1, x_2} \in A$ be any group element such that $\tilde{b}_{x_1, x_2} \cdot P_{1} = P_{2}$. Finally, let $b_{x_1, x_2}(P \to C)$ be the coset of $\tilde{b}_{x_1, x_2}$.
\end{defn}
The following lemma follows immediately from \cite[Proposition 2.4]{BCK}.
\begin{lemma}\label{transition-data}
Let $(P \to C), (P' \to C') \in \Mbar_{0, X }^A(\ell_X)$. Then the following statements are equivalent.
\begin{enumerate}
\item The points $(P \to C)$ and $(P' \to C')$ are in the same connected component of $\Mbar_{0, X }^A(\ell_X)$.
\item For all $x_1, x_2 \in X$, 
\[ b_{x_1, x_2}(P \to C) = b_{x_1, x_2}(P' \to C'). \]
\item For a fixed element $z \in X$,
\[ b_{z, x}(P \to C) = b_{z, x}(P' \to C') \]
for all $x \in X$.
\end{enumerate}
\end{lemma}

 Given a function $\ell_X: X \to A$ and an element $a \in A$ such that
\[ \sum_{x \in X} \ell_X(x) + a = 0, \]
we write
\[\Mbar_{0, X \sqcup \{*\}}^A(\ell_X, -a) \]
for the moduli space of genus zero $X \sqcup \{*\}$-marked admissible $A$ covers, with monodromy at the marked point $p_x$ given by $\ell_X(x)$ and monodromy at $p_*$ given by $-a$. Lemma \ref{transition-data} motivates the following definition. 
\begin{defn}
Given a point $(P \to C) \in \Mbar_{0, X \sqcup \{* \}}^A(\ell_X, -a)$, a tuple of elements
\[ (b_{*, x} \in A \mid x \in X) \]
satisfying
\[ b_{*, x} \cdot \langle \ell_X(x), a \mid x \in X \rangle = b_{*, x}(P \to C)\] is
a \textit{\textbf{transition datum}} for the point $P \to C$. 
\end{defn}
As per Lemma \ref{transition-data}, a transition datum for a cover $(P \to C)$ specifies the connected component of $\Mbar_{0, X \sqcup \{*\} }(\ell_X, -a)$ which contains it.

Now we explain the action of $\tFSA\op$ on the homology groups $H_*(\Mbar_{g, X}^A(\ell_X))$. Given a morphism
\[(f, g): (X, \ell_X) \to (Y, \ell_Y) \]
in $\tFSA\op$, we get, for each $y \in Y$, a connected component $\mathcal{T}_y$ of the moduli space
\[\Mbar_{0, f^{-1}(y) \sqcup y}^A(\ell_{X}|_{f^{-1}(y)}, -\ell_Y(y)), \]
by choosing the transition datum determined by the function \[g|_{f^{-1}(y)}: f^{-1}(y) \to A; \]
i.e. we have
\[b_{y, x} = g(x) \]
for each $x \in f^{-1}(y)$.
Choose for each $y \in Y$ a pointed admissible $A$-cover \[(P_y, (\tilde{s}_{q})_{q \in f^{-1}(y)}, \tilde{r}_y) \to (C_y, (s_q)_{q \in f^{-1}(y)}, r_y) \in \mathcal{T}_y. \]
Given a pointed admissible $A$-cover
\[ (Q, (\tilde{p}_y)_{y \in Y}) \to (D, (p_y)_{y \in Y})) \in \Mbar_{g, Y}^{A}(\ell_Y),  \]
we obtain a pointed admissible $A$-cover
\[ \left( \left(Q  \sqcup \coprod_{y \in Y} P_y\right)/\sim, (\tilde{s}_x)_{x \in X} \to \left(D \sqcup  \coprod_{y \in Y} C_y \right)/\sim, (s_x)_{x \in X} \right)  \in \Mbar_{g, X}^A(\ell_X) , \]
where the equivalence relations are given by gluing the points $\tilde{p}_y$ to $\tilde{r}_y$ in the source, and gluing $p_y$ to $r_y$ in the target. In this way we have defined, given our choices of points in the connected components $\mathcal{T}_y$, a morphism
\[ \varphi_{(f, g)^*} : \Mbar_{g, Y}^{A}(\ell_Y) \to \Mbar_{g, X}^{A}(\ell_X), \]
which in turn defines a morphism
\begin{equation}
(f, g)^*: H_i(\Mbar_{g, Y}^{A}(\ell_Y)) \to H_i(\Mbar_{g, X}^{A}(\ell_X)) 
\end{equation}
on the level of homology. Since the connected components $\mathcal{T}_y$ are uniquely determined by $(f, g)$, the morphism $(f, g)^*$ is well-defined. Now we must check that the association $(f, g) \mapsto (f, g)^*$ defines an action of the category $\tFSA\op$.
\begin{prop}
Let $(f_1, g_1): (X, \ell_X) \to (Y, \ell_Y)$ and $(f_2, g_2):(Y, \ell_Y) \to (Z, \ell_Z)$ be morphisms in the category $\tFSA\op$. Then
\[ \left((f_2, g_2) \circ (f_1, g_1)\right)^* = (f_1, g_1)^* \circ (f_2, g_2)^* \]
\end{prop}
\begin{proof}
We want to verify that the two morphisms determine the same connected component of the moduli space
\[ \prod_{z \in Z} \Mbar_{0, (f_2 \circ f_1)^{-1}(z) \sqcup \{z\}}(\ell_X, -\ell_Z(z)). \]
It suffices to verify that they determine the same connected component of
\[ \Mbar_{0, (f_2 \circ f_1)^{-1}(z) \sqcup \{z\}}(\ell_X, -\ell_Z(z)) \]
for each $z \in Z$. By Lemma \ref{transition-data}, it suffices to show that we can find transition data \[(b_{z, x} \in A \mid x \in (f_2 \circ f_1)^{-1}(z) ) \]
for the connected component determined by $((f_2, g_2) \circ (f_1, g_1))^*$ and
\[(c_{z, x} \in A \mid x \in (f_2 \circ f_1)^{-1}(z)) \]
for the connected component determined by $(f_2, g_2)^* \circ (f_1, g_1)^*$ such that $b_{z, x} = c_{z, x}$ for all $x \in (f_2 \circ f_1)^{-1}(z)$. On one hand, the morphism $\left((f_2, g_2) \circ (f_1, g_1)\right)^*$ has transition datum
\[ b_{z, x} = g_1(x) + (g_2 \circ f_1)(x). \]
We want to verify that the same transition datum specifies the connected component obtained by acting first by $(f_2, g_2)^*$ and then by $(f_1, g_1)^*$. First, $(f_2, g_2)$ specifies the connected component $\mathcal{T}_z$ of
\[ \Mbar_{0, f_2^{-1}(z) \sqcup \{z\}} (\ell_Y, -\ell_Z(z)) \]
with transition datum $b_{z, y} = g_2(y)$ for each $y \in f_2^{-1}(z)$. Then, for each $y \in Y$, $(f_1, g_1)$ specifies the connected component $\mathcal{T}_{y}$ of
\[ \Mbar_{0, f_1^{-1}(y) \sqcup \{y\}}^{A}(\ell_X, -\ell_Y(y))  \]
with transition datum $b_{y, x} = g_1(x)$ for each $x \in f_1^{-1}(y)$. Now fix $y \in f_2^{-1}(z)$, and for a cover $P \to C$ underlying a point in $\mathcal{T}_z$
write $\hat{p_y}$ for the chosen lift of $p_y$ to $P$, and $\hat{p_z}$ for the chosen lift of $p_z$ to $P$. Given a cover $P' \to C'$ underlying a point in $\mathcal{T}_y$,
write $\tilde{p_y}$ for the chosen lift of $p_y$ to $P'$. Fix $x \in f_1^{-1}(y)$, and write $\tilde{p_x}$ for the chosen lift of $p_x$ to $P'$. The connected component specified by $(f_1, g_1)^* \circ (f_2, g_2)^*$ is given by gluing $C$ to $C'$ at their respective marked points $p_y$, and gluing their fibers in $P$ and $P'$ such that $\hat{p_y}$ is glued to $\tilde{p_y}$. Let $\overline{b}$ denote the transition data of this glued cover; we will compute $\overline{b}_{z, x}$. Since we are gluing $\hat{p_y}$ to $\tilde{p_y}$, we have the following: if $\alpha \in A$ is a group element whose action takes the connected component of $P$ containing $\hat{p_z}$ to the connected component containing $\hat{p_y}$, and $\beta \in A$ is a group element whose action takes the connected component of $P'$ containing $\tilde{p}_y$ to the connected component containing $\tilde{p}_x$, then the group element $\alpha + \beta \in A$ takes the connected component of the glued curve $P \sqcup P'$ containing $\hat{p}_z$ to the connected component containing $\tilde{p}_x$. Since we can take $\alpha = g_2(y)$ and $\beta = g_1(x)$, we have that
\[\overline{b}_{z, x} = g_1(x) + g_2(y) = g_1(x) + (g_2 \circ f_1)(x),  \]
which completes the proof.
\end{proof}

\section{Finite generation}\label{sec:finite_generation}

Now that we have proven the assignment, 
\[\cV_{i,g}: (X, \ell_X) \mapsto H_i(\overline \cM_{g,X}^A(\ell_X); \bbQ)\] is an $\tFSA\op$-module, we may prove it has finite generation degree as an $\tFSA\op$-module.

\subsection{The generation degree of $\mathcal{V}_{0,0}$} The starting point for our proof of finite generation is expanding the definition of $\mathcal{V}_{0,0}$ and showing it is finitely generated in degree 1. Since the space $\overline{\mathcal{M}}_{0, X}^A(\ell_X)$ is only defined when $|X| \geq 3$, we need also give a definition of $\cV_{i, 0}$ when $|X| < 3$. To do so, we use the classification of the connected components of $\overline{M}_{0, X}^A(\ell_X)$.

\begin{defn}
    Let $\cV_{0,0}$  be the free $\tFSA\op$-module generated (as a vector space) by the $\tFSA\op$-set  $\mathcal{S}_{0, 0}$. $\mathcal{S}_{0,0}$ is defined on objects by
\[\mathcal{S}_{0, 0}(X, \ell_X) := \left(A / \langle \ell_X(x) \mid x \in X \rangle\right)^X / A.  \]
\end{defn}

Note that with this definition, $\cV_{0,0}$ coincides with the existing definition of $\cV_{i,g}: (X, \ell_X) \mapsto H_i(\overline \cM_{g,X}^A(\ell_X); \bbQ)$ when $|X| \geq 3$. To describe the $\tFSA\op$-module structure of $\mathcal{V}_{0, 0}$, it suffices to describe the how morphisms act on $\mathcal{S}_{0, 0}$.
Suppose we are given a morphism $(f, g): (X, \ell_X) \to (Y, \ell_Y)$ in $\tFSA$. Similarly to the definition of the $\tFSA\op$-module $\tV$, $f$ defines a function
\[h: \left(A / \langle \ell_X(y) \mid y \in Y\rangle\right)^Y \to \left(A / \langle \ell_X(x) \mid x \in X \rangle\right)^X,  \]
which can be checked to be compatible with the action of $A$. Finally, $h$ is perturbed by the pointing $g$: given a function $\alpha \in \left(A / \langle \ell_Y(y) \mid y \in Y\rangle\right)^Y$, we define the function \[(f, g)^*(\alpha): X \to \frac{A}{\langle \ell_X(x) \mid x \in X\rangle}\]
by
\[ h(\alpha)(x) + g(x) \]
Given this combinatorial description of $\mathcal{V}_{0, 0}$, we now prove that it is generated in degree 1.

\begin{thm}\label{v00}
The $\tFSA\op$-set $\mathcal{S}_{0, 0}$, and hence the $\tFSA\op$-module $\mathcal{V}_{0, 0}$, is finitely generated in degree $1$. 
\end{thm}

\begin{proof}
The objects of $\tFSA$ in degree 1 are pairs $(\{*\}, a)$, where $\{*\}$ is a set with one element and $a \in A$ is the labeling of this element. Given a general object $(X, \ell_X)$ of $\tFSA$, there exists a unique element $a \in A$ such that $(X, \ell_X)$ admits a morphism to $(\{*\}, a)$, determined by $a = \sum_{x \in X} \ell_X(x)$. Such a morphism is of the form $(f, g)$, where $f : X \to \{*\}$ is the constant map and $g: X \to A$ is any pointing. Let us compute $(f, g)^*(\alpha)$ where $\alpha$ is the unique element of
\[ \mathcal{S}_{0, 0}(\{*\}, a) = (A / \langle a \rangle)^{\{*\}}/A . \] First, choose a function $c: \{*\} \to A/\langle a \rangle$ lifting the function $\alpha$: all choices $\{*\} \to A$ are possible, and all choices are rendered equivalent after quotienting by the action of $A$. Precomposing with $f$ and postcomposing by the quotient map $A /\langle a \rangle \to A/\langle \ell_X(x) \mid x \in X\rangle$ gives us a constant function $h : X \to A/\langle \ell_X(x) \mid x \in X \rangle$. Finally, $h$ is perturbed by $g$, giving
\[(f, g)^*(\alpha)(x) = g(x) + h(x). \]
Since $h(x)$ is constant, by varying our choice of the pointing $g$, it is possible to express \textit{any} function $X \to A/ \langle \ell_X(x) \mid x \in X \rangle$ as $(f, g)^*(\alpha)$ for some $g$. Hence, $\mathcal{S}_{0, 0}$ is finitely generated in degree 1.
\end{proof}

Let $(P \to C) \in \Mbar_{g,X}^A(\ell_X)$ and  $(Q \to D) \in \Mbar_{h,Y}^A(\ell_Y)$ be pointed admissible $A$-covers.  Given $x \in X$ and $y \in Y$ such that $\ell_X(x) + \ell_Y(y) = 0$, we may form the \emph{glued cover}, $${\text{glue}}_{x,y}(P \to C,Q \to D) \in \Mbar_{g+h, X - x \sqcup Y - y} ^A(\ell_X|_{X - x} \sqcup \ell_Y|_{Y - y}),$$ obtained by gluing $(P \to C)$ and $(Q \to D)$ along the unique isomorphism of $A$-bundles between $P|_x \to x$ and $Q|_y \to y$ taking the lift of $x$ to the lift of $y$. In particular, for any $a \in A$ there is a gluing map 
\begin{equation}
   \text{glue}_a:  \Mbar_{g, X \sqcup *} (\ell_X \sqcup a)  \times \Mbar_{g, Y\sqcup *} (\ell_X \sqcup -a)  \to \Mbar_{g+h, X \sqcup Y} (\ell_X \sqcup \ell_Y).\label{gluing}
\end{equation}

Lemma \ref{transition-data} implies that there is an isomorphism of $\tFSA\op$-modules
$\mathcal{V}_{0, 0} \cong \overline{V}_0$. Using this isomorphism, we obtain a quotient map $q: \tV_0 \to \mathcal{V}_{0,0}$. A key step in our inductive argument for finite generation of the modules $\mathcal{V}_{i, g}$ is the following lemma, as it relates the natural suspension map of Proposition \ref{SuspensionMapProp} to this gluing map via the quotient map $q$. 

\begin{lemma}\label{factorgluing}
The natural map
\[\tilde{\eta}: \bigoplus_{a \in A} \Sigma_a \tV_{0} \oast \Sigma_{-a} \mathcal{V}_{i, g} \to \mathcal{V}_{i, g}\]
factors as
\[
\begin{tikzcd}
\bigoplus_{a \in A} \Sigma_a \tV_{0} \oast \Sigma_{-a} \mathcal{V}_{i, g} \arrow[rr] \arrow[d, "q"]  & & \mathcal{V}_{i, g} \\
\bigoplus_{a \in A} \Sigma_a \mathcal{V}_{0,0} \oast \Sigma_{-a} \mathcal{V}_{i, g} \arrow[urr, swap, "\mathrm{glue}"]
\end{tikzcd}.
\]
\end{lemma}
\begin{proof}
By definition,
\[\bigoplus_{a \in A} \Sigma_a \tV_{0} \oast \Sigma_{-a} \mathcal{V}_{i, g} = \bigoplus_{a \in A} \bigoplus_{X = S \sqcup T} \Sigma_{a} \tV_{0}(S, \ell_{X}|_{S}) \otimes \Sigma_{-a}\mathcal{V}_{i, g}(T, \ell_{X}|_{T}). \]
It suffices to prove that the diagram commutes for some fixed $a \in A$ and partition $X = S \sqcup T$; that is, we would like to show that the diagram
\begin{equation}\label{qfactors}
\begin{tikzcd}
 \Sigma_a \tV_{0}(S, \ell_X|_S) \otimes \Sigma_{-a} \mathcal{V}_{i, g}(T, \ell_X|_T) \arrow[rr] \arrow[d, "q"]  & & \mathcal{V}_{i, g}(X, \ell_X) \\
\Sigma_a \mathcal{V}_{0,0}(S, \ell_X|_S) \otimes \Sigma_{-a} \mathcal{V}_{i, g}(T, \ell_X|_T) \arrow[urr, swap, "\mathrm{glue}"]
\end{tikzcd}
\end{equation}
commutes, for any $a \in A$ and partition $X = S \sqcup T$. We may assume that $|S|, |T| \geq 2$, and that $\sum_{x \in S} \ell_X(s) = -a$, so that the vector spaces on the left-hand side of the diagram are nonzero. We may rewrite 
\begin{align*}
\Sigma_a \tV_{0}(S, \ell_X|_S) &= \mathbb{Q} \cdot A^{S \sqcup \{*\}}/A \\
\Sigma_a \mathcal{V}_{0, 0}(S, \ell_X|_S) &=  \mathbb{Q} \cdot \left(A/\langle a, \ell_X(s) \mid s\in S\rangle\right)^{S \sqcup \{*\}}/A,
\end{align*}
so the map $q$ on the left-hand side is induced by the quotient map
\[A^{S \sqcup \{*\}}/A \to  \left(A/\langle a, \ell_X(s) \mid s\in S\rangle\right)^{S \sqcup \{*\}}/A. \] 

Then the diagram in (\ref{qfactors}) becomes 

\[
\begin{tikzcd}
 \bigoplus_{g: S \to A} \mathcal{V}_{i, g}(T \sqcup \{*\}, \ell_X|_T \sqcup -a) \arrow[rr, "{\oplus_{g: S \to A}(f_S, \overline{g})^*}" ] \arrow[d, "q"]  & & \mathcal{V}_{i, g}(X, \ell_X) \\
\bigoplus_{c \in \pi_0(\Mbar_{0, S \sqcup \{*\}}^A (\ell_X|_{S}, a) ) } \mathcal{V}_{i, g}(T \sqcup \{*\}, \ell_X|_T \sqcup -a) \arrow[urr, swap, " {\oplus_{c} (\mathrm{glue}_c})_*"]
\end{tikzcd}
\]

with $\tFSA\op$ morphism ($f_S: X \to T \sqcup \{*\}, \overline{g}: X \to A$) given by

\begin{equation*}
\begin{aligned}[c]
f_S(x) = \begin{cases} * \text{ if } x \in S\\ x \text{ if } x \in X \setminus S \end{cases}
\end{aligned}
\qquad\qquad
\begin{aligned}[c]
\overline{g}(x) = \begin{cases}g(x) \text{ if } x \in S\\ 0 \text{ if } x \in X \setminus S \end{cases}
\end{aligned}
\end{equation*}

To check commutativity, fix a pointing $g: S \to A$. This defines a single summand in the top left corner of the diagram. For any $v \in \mathcal{V}_{i, g}(T \sqcup \{*\}, \ell_X|_T \sqcup -a)$, the horizontal map applies $(f_S, \overline{g})^*$ to $v$. By the definition of the $\tFSA\op$ structure on $\mathcal{V}_{i, g}$, this means that we first find the connected component \[c \in \pi_0\left(\Mbar_{0, S \sqcup \{*\}}^A (\ell_X|_{S}, a)\right)\]
determined by $(f_S, \overline{g})$, choose a point inside of this connected component, and then apply to $v$ the linearization of the gluing map induced by gluing the curve represented by the point. The described procedure is precisely $(\mathrm{glue}_c)_* \circ q$ for the specified component $c$. 
\end{proof}

 Using the gluing map defined in \eqref{gluing}, one obtains a map $$\glue_a(i_1,i_2,g_1,g_2):\Sigma_{a} \cV_{i_1,g_1} \oast \Sigma_{-a} \cV_{i_2,g_2} \to \cV_{i_1 + i_2, g_1 + g_2}.$$ We then consider the following sum of such gluing maps, 

 \begin{equation}\label{tau def}
    \tau_{i,g}:\bigoplus_{a \in A}  \left(\Sigma_a(\Sigma_{-a} \cV_{i, g-1}) \oplus \bigoplus_{g_1 + g_2 = g} \bigoplus_{i_1 + i_2 = i} \Sigma_{a} \cV_{i_1,g_1} \oast \Sigma_{-a} \cV_{i_2,g_2}\right) \longrightarrow \cV_{i, g}
\end{equation}

 Defined as
 
 \[\tau_{i,g} = \bigoplus_{a \in A} \left(\glue_a \oplus \bigoplus_{i_1 + i_2 = i} \bigoplus_{g_1 + g_2 = g} \glue_a(i_1, g_1, i_2, g_2)\right).\]

 Once applied to an object, these maps can be interpreted geometrically as gluing maps of admissible $A$-covers. For an $A$-cover, summands $\Sigma_a \Sigma_{-a} \tV_{i, g-1}$ correspond to gluing two points in the same connected component of the target curve into a node, thus increasing the genus by 1. The other summands correspond to gluing together two $A$-covers as previously described.

\begin{prop}\label{surjective}
    The map $\tau_{i,g}$ is surjective in all degrees at least $i + 4$ if $g = 0$, and it is surjective in all degrees at least $\max(i - 2g + 3, 1)$ if $g \geq 1$.
\end{prop}
\begin{proof}
Let $\partial \overline \cM_{g,X}^A( \ell_X)$ denote the boundary divisor to the inclusion  $\cM_{g,X}^A( \ell_X) \subseteq \overline \cM_{g,X}^A(\ell_X)$.   The long exact sequence in Borel--Moore homology gives an exact sequence 
$$H_i(\partial \overline \cM_{g,X}^A( \ell_X)) \to H_i( \overline \cM_{g,X}^A( \ell_X)) \to H_i^{\rm BM}(\cM_{g,X}^A( \ell_X) ).$$ Applying Proposition \ref{vanish} gives that the first map is surjective when $|X| \geq i + 4$ if $g = 0$ and when $|X| \geq \max(i - 2g + 3,1)$ if $g \geq 1$.  Because $H_i( \overline \cM_{g,X}^A(\ell_X))$ is pure of weight $2i$, it  follows that the pure component of $H_i(\partial \overline \cM_{g,X}^A(\ell_X))$ surjects onto $H_i( \overline \cM_{g,X}^A(\ell_X))$. Now the map $\tau_{g,i}$ is induced on homology by the disjoint union of the gluing maps \[
\begin{tikzcd}
\coprod_{a \in A} \left(\Mbar_{g - 1, X \sqcup \{\star, \bullet \}}^A(\ell_{X} \sqcup \{a, -a\}) \sqcup \coprod_{ X_1 \sqcup X_2 = X} \Mbar_{g_1, X_1 \coprod \{\star \}}^A ({\ell_X}|_{X_1} , a) \times \Mbar_{g_2, X_2 \sqcup \{\bullet\}}^A({\ell_X }|_{X_2} , - a)\right) \arrow[d] \\ \partial\Mbar_{g, X}^A(\ell_X),
\end{tikzcd}\]
followed by the inclusion of the boundary. By definition, the above map surjects onto $\partial \overline \cM_{g,X}^A(\ell_X)$, hence the induced map on the pure part of the homology surjects onto the pure component of $H_i(\partial \overline \cM_{g,X}^A(X, \ell_X))$ by   
\cite[Lemma A.4]{lewis1999survey}. Therefore, when $g = 0$ and $|X| \geq i + 4$, or when $g \geq 1$ and $|X| \geq \max(i-2g +3, 1)$, the composite surjects onto $H_i(\overline \cM_{g,X}^A(\ell_X))$.
\end{proof}

To obtain a bound on when $\tau_{i,g}$ is surjective that is not dependent on whether $g = 0$, we prove the following numerical lemma. 

\begin{lemma}\label{bound}
    Let $f(i,g): \bbN^2 \to \bbN$ be a function satisfying $f(0,0) = 3$, $f(0,1) = 1$ and $f(0,g) = 0$ for $g \geq 2$.   Suppose that \[f(i,0) \leq \max\left(i + 4,  ~\max_{\substack{(i_1, g_1) + (i_2,g_2) = (i,g) \\ (i_k,g_k) \neq (0,0)}} (f(i_1, g_1) + f(i_2, g_2))\right) \] 
    for $i > 0$ and that
    \[f(i,g) \leq \max\left(\max(i -2g +3, 1),  f(i, g- 1), ~\max_{\substack{(i_1, g_1) + (i_2,g_2) = (i,g) \\ (i_k,g_k) \neq (0,0)}} (f(i_1, g_1) + f(i_2, g_2))\right) \]
    if $g, i \geq 1$. Then $f(i,g) \leq g + 5i$ if $(i, g) \neq (0, 0)$.
\end{lemma}
\begin{proof}
    We will prove by induction that $f(i, g) \leq g + 5i$ when $(i, g) \neq (0, 0)$; the statement of the lemma follows. The base cases will be $(0, 1)$ and $(1, 0)$. The case of $(0, 1)$ is clear, and for $(i, g) = (1, 0)$ we have that
    \[ i + 4 = 5i, \]
    proving the claim. Now suppose the claim holds for all $(i', g')$ which are lexicographically smaller than $(i, g)$, and suppose that $g > 0$. We can also suppose that $i > 0$, since the claim is direct when $i = 0$. In this case, $i - 2g + 3 \leq g + 5i$ since $g > 0$, and $1 \leq g + 5i$, and $f(i, g - 1) \leq g - 1 + 5i \leq g+ 5i$, by induction. Also by induction,
    \[\max_{\substack{(i_1, g_1) + (i_2,g_2) = (i,g) \\ (i_k,g_k) \neq (0,0)}} (f(i_1, g_1) + f(i_2, g_2)) \leq \max_{\substack{(i_1, g_1) + (i_2,g_2) = (i,g) \\ (i_k,g_k) \neq (0,0)}} g_1 + 5i_1 + g_2 + 5i_2 = g + 5i, \]
    proving the claim. The case $g = 0$ is similar.
\end{proof} 

We now give an inductive argument on $\cV_{i,g}$ to prove its finite generation. 

\begin{thm} \label{maintheorem}
    If $(i, g) \neq (0,0)$, the $\tFSA \op$-module $\cV_{i,g}$ is finitely generated in degree $\leq \mainbound$.
\end{thm}
\begin{proof}

     We induct lexicographically on the tuple $(i, g)$. Our proof proceeds by constructing a quotient $\tFSA\op$-module $Q_{i,g}$ of $\cV_{i,g}$ whose finite generation implies finite generation of $\cV_{i,g}$. 
     
    Recall the gluing map
    
    \[\tau_{i,g}:\bigoplus_{a \in A}  \left(\Sigma_a(\Sigma_{-a} \cV_{i, g-1}) \oplus \bigoplus_{g_1 + g_2 = g} \bigoplus_{i_1 + i_2 = i} \Sigma_{a} \cV_{i_1,g_1} \oast \Sigma_{-a} \cV_{i_2,g_2}\right) \longrightarrow \cV_{i, g}.\]
    Applying Proposition \ref{surjective} and Lemma \ref{bound} implies that $\tau_{i,g}$ is surjective when applied to $(X, \ell_X)$ such that $|X| > g + 5i$. Our goal is to leverage surjectivity of $\tau_{i,g}$ to imply finite generation of $\cV_{i,g}$ from finite generation of the domain.

    To define the quotient $Q_{i,g}$, consider the summands $\Sigma_a \cV_{i_1, g_1}\oast \Sigma_{-a}\cV_{i_2, g_2}$ of the domain of $\tau_{i,g}$ such that neither $(i_1, g_1)$ nor $(i_2, g_2)$ is identically $(0,0)$. In the case when $(i_k, g_k) = (0,0)$, the corresponding summands in the left-hand side are of the form $\Sigma_a \cV_{0, 0} \oast \Sigma_{-a}\cV_{i,g}$ or $\Sigma_a\cV_{i,g} \oast \Sigma_{-a} \cV_{0,0}$. The images of these summands will be all of $\cV_{i,g}$, so in the cases $(i,g) = (0,1)$ or $(1,0)$, $\tau_{i,g}$ is immediately surjective. 
    
    We may define a sub-$\tFSA \op$-module, $I_{i,g}$, by taking the image of these summands under $\tau_{i,g}$:

    \begin{equation}\label{Idef}
        I_{i,g} := \text{im}_{\tau_{i,g}} \left(\bigoplus_{a \in A} \left(\Sigma_{a, -a}\cV_{i, g-1} \oplus\bigoplus_{\substack{(i_1, g_1) + (i_2, g_2) = (g,i) \\ (i_k, g_k) \neq (0,0)}} \Sigma_a \cV_{i_1, g_1} \oast \Sigma_a\cV_{i_2, g_2} \right)\right).
    \end{equation}
    Then let \[Q_{i,g} := \frac{\cV_{i,g}}{I_{i,g}}.\]

    For an $\tFSA \op$-module $M$, we will write $\deg(M)$ for the generation degree of $M$. Since $Q_{i,g}$ is a quotient module, and shifting does not affect generation degree, Propositions \ref{restricted} and \ref{quotient} give immediately that  $\deg(\cV_{i,g}) \leq \max(\deg(Q_{i,g}), \deg(I_{i,g})).$ Furthermore, applying Proposition \ref{product} to each summand in $I_{i,g}$, we have that \[\deg(I_{i,g}) \leq \max_{\substack{(i_1, g_1) + (i_2, g_2) = (i,g) \\ (i_k, g_k) \neq (0,0)}} (\deg(\cV_{i_1, g_1}) + \deg(\cV_{i_2, g_2})).\] 

    We now proceed by induction to prove surjectivity of $\tau_{i,g}$ in the desired degrees. The base cases of our induction are $(i, g) = (0, 0), (0,1)$, and $(1,0)$. By Theorem \ref{v00}, $(0,0)$ is already finitely generated in degree 1, and as noted previously, $\tau_{i,g}$ is immediately surjective if $(i,g) = (0,1), (1,0)$. 
    
    Since $I_{i,g}$ is the image of only $\cV_{i_k, g_k}$ such that $i_k + g_k > 0$, by induction both  $\deg(\cV_{i_1, g_1})$ and $\deg(\cV_{i_2, g_2}) \leq g + 5i$. Thus, $\deg(I_{i,g}) \leq g + 5i$. What remains is to prove that $\deg(Q_{i,g}) \leq g + 5i$ as well. 

    Recall that for any $\tFSA\op$-module $M$, there is a natural map of shifted modules 
    \[\tilde{\eta}_M: \bigoplus_{a \in A} \Sigma_{-a}\tV_{0} \oast \Sigma_a M \to M.\]
    By Proposition \ref{SuspensionMapProp}, finite generation in degree $\leq d$ of $M$ may be equivalently characterized by $\tilde{\eta}_M$ being surjective when applied to any set of size $> d$. We may define identically a map $\tilde{\eta}'_{M}$: 

    \[\tilde{\eta}'_{M}: \bigoplus_{a 
    \in A} \Sigma_a M \oast \Sigma_{-a}\tV_0 \to M\]
    with the same definition as in (\ref{eta tilde}), just swapping the labels $a$ and $-a$. Letting $\rho$ denote the induced map from the quotient $q: \cV_{i,g} \to Q_{i,g}$, we obtain the following diagram of $\tFSA\op$-modules.

    \begin{equation}\label{Qdiagram}
    \begin{tikzcd}
        \bigoplus_{a \in A} \Sigma_a \tV_0 \oast \Sigma_{-a} \cV_{i, g} \oplus \bigoplus_{a \in A}\Sigma_{a}\cV_{i,g} \oast \Sigma_{-a} \tV_0  \arrow[rr, "\tilde\eta_{\cV_{i,g}} \oplus \tilde\eta_{\cV_{i,g}}'"] \arrow[d, "\rho"]& &\cV_{i,g} \arrow[d, "q"]\\
        \bigoplus_{a \in A} \Sigma_a \tV_0 \oast \Sigma_{-a} Q_{i,g} \oplus \bigoplus_{a \in A}\Sigma_{a} Q_{i,g} \oast \Sigma_{-a} \tV_0 \arrow[rr, "\tilde{\eta}_{Q_{i,g}} \oplus \tilde{\eta}_{Q_{i,g}}'"] & &Q_{i,g}
    \end{tikzcd}
    \end{equation}
    Varying $M$ over all $\tFSA\op$-modules, the maps $\tilde{\eta}_M$ define the component maps of a natural transformation from the functor $M \mapsto \oplus_{a \in A}\Sigma_{-a} \tV_0 \oast \Sigma_a M$ to the identity functor. (Similarly for $\tilde{\eta}'$.) Commutativity of (\ref{Qdiagram}) is then immediate from this naturality.

    The image of $\tilde \eta_{\cV_{i,g}} $ equals the image of $\tilde \eta_{\cV_{i,g}}'$, essentially because the image of the linearization of a map of spaces $f: X \times Y \to Z$ is the same as the image after composing with $Y \times X \to X \times Y$. It follows that the image of $\tilde{\eta}_{Q_{i, g}}$ coincides with that of $\tilde{\eta}_{Q_{i, g}}'$, and that the cokernel of this map is the same as $\mathrm{coker}(\tau_{i,g})$. It follows that $\tilde{\eta}_{Q_{i, g}}$ is surjective whenever $\tau_{i,g}$ is, and in particular when $|X| > g + 5i$ by Proposition \ref{surjective}. Thus by Proposition \ref{SuspensionMapProp}, $Q_{i, g}$ is generated in degree $\leq g+ 5i$, and the proof is complete.
    \end{proof}

    \section{Further questions}\label{sec:further questions}
We list some questions motivated by our work.

\subsection{Non-abelian generalizations} The definition of the categories $\FSA$ and $\tFSA$ only makes sense for $A$ an abelian group. Are there analogs of these categories for general finite groups that can be used to study $\overline \cM_{g,n}^G$? Since the category $\FSA$ is related to the components of $\overline \cM_{0,n}^A$, a natural first step is to classify the connected components of $\Mbar_{0, n}^G$ for arbitrary $G$, and attempt to build a category from these components.

\subsection{Refined results on Hilbert series} Given a finitely generated $\tFSA\op$ module $M$, and any $\ell: X \to A$ where $X$ is a finite set, the vector space $M(X, \ell)$ is naturally a representation of the  semidirect product ${\rm Aut}(X,\ell) \ltimes A^X$, where ${\rm Aut}(X,\ell)$ is the subgroup of permutations of $X$ that fix $\ell$.  In the case where $A = e$, we obtain a sequence of $S_n$ representations, and the form of the associated symmetric function was determined in \cite{tosteson2022algebra}.   Is there an analogous theory for general abelian groups?   Such a theory would lead to enhancements of Theorem A and Corollary D.

\subsection{Explicit computations}  We establish the qualitative form of the generating functions for the homology of spaces of admissible covers.   It would be interesting to compute these generating functions, at least for small values of $g,i$.  

A related problem is to find explicit presentations, by generators and relations, of the $\FSA \op$  and $\tFSA \op$ modules that we construct.  We bound the generation degree of these modules, is there a similar bound on presentation degree?  If so does it bound the degree of the numerators in Theorem A?

    \bibliographystyle{alpha}
\bibliography{AdmBibliography}

\end{document}